\documentstyle[11pt]{article}
\renewcommand{\baselinestretch}{1.3}

\newtheorem{conj}{conjecture}

\def\TB{{\bf b^\infty}}
\def\TM{{\bf m^\infty}}

\def\deq{\, {\stackrel {def} {=}}}
\def\diseq{\, {\stackrel {{\cal D}} {=}}}

\def\dd{\delta}
\def\ee{\epsilon}
\def\E{{\bf{E}}}
\def\P{{\bf{P}}}
\def\N{\hbox{I\kern-.2em\hbox{N}}}
\def\R{\hbox{I\kern-.2em\hbox{R}}}
\def\AB{{\overline{A}}}

\def\C{{\cal{C}}}
\def\Z{{\bf{Z}}}

\def\F{{\cal{F}}}
\def\cE{{\cal{E}}}

\def\cEB{\overline{\cal{E}}}
\def\|{\, | \, }

\def\Capb{\overline{{ \rm cap}}}
\def\Cap{{\rm cap}}

\begin{document}
\begin{titlepage}
\begin{center}
{\large \bf TREE-INDEXED PROCESSES} \\
\end{center}
\vspace{5ex}
\begin{flushright}
Robin Pemantle \footnote{Department of Mathematics, University of 
Wisconsin-Madison, Van Vleck Hall, 480 Lincoln Drive, Madison, WI 53706}$^,$
\footnote{Research supported in part by National Science Foundation Grant 
\# DMS 9300191, by a Sloan Foundation Fellowship and by a Presidential
Faculty Fellowship} ~\\
\end{flushright}

\vfill

ABSTRACT:

\noindent{This} article examines a recent body of work on 
stochastic processes indexed by a tree.  Emphasis is on the
application of this new framework to existing probability models.
Proofs are largely omitted, with references provided.  

\vfill



\end{titlepage}

\section{Introduction}

Tree-indexed processes are not really new stochastic processes, but
rather new ways of looking at already existing probability models.
Consider, by way of analogy, ordinary continuous-time stochastic processes.  
These are of course ``merely'' collections of random variables indexed
by the positive reals.  But when viewed as random trajectories, powerful
concepts such as filtrations and stopping times naturally arise, which
are fruitful -- indeed necessary -- for successful analysis
of the original problems.  Similarly, many probability models
involving trees may be described as follows.  First pick a tree,
either deterministically or at random.  Then attach some randomness
to the tree (think of real random variables on each edge or vertex) and 
ask questions about the resulting structure.  The tree-indexed viewpoint
is to think of this as a random field indexed by the space of
paths through the tree and taking values in the space of sequences
of real numbers\footnote{As far as I know this viewpoint dates from
1990 when the preprint of Evans (1992) was circulated.}.

My main concern in this article is to illustrate how the tree-indexed
view may be applied to a variety of well-known models, and to show
how some general theory may be used to extract information about
these models in a relatively painless way.  Just as potential theory 
(the study of potential, energies and capacities) is
almost synonymous with the classical theory of Markov processes,
the potential theory of trees is behind most of the theorems surveyed
here, and in fact I will not draw a distinction between tree-indexed
theory and potential theory or second moment methods on trees.  Since 
trees are easier to analyze than lattices, there are many papers proving
results on trees as a somewhat unmotivated alternative or
a ``high-dimensional analogue'' to Euclidean space\footnote{\em Mea culpa}.  
In this survey I will emphasize models where the trees are there
because nature put them there.  I will also discuss several applications
of tree-indexed processes to questions that do not appear at first glance
to involve trees.  Indeed the applications of tree-indexed theory to
the intersections of random subsets of Euclidean space via the 
tree-representation of $[0,1]^d$ are some of the most compelling 
justifications of tree-potential theory.

There are two ways I can indicate the scope of this survey.  One is to
begin by listing the models and the questions that are addressed by
tree-indexed theory.  The other is to state the basic definitions
and the fundamental theorems.  This section takes the
former approach, discussing questions susceptible to tree-indexed
theory.  These questions predate by far the emergence of tree-indexed
theory, so much of the background given here is quite classical.  
In particular, problems in branching process theory and fractal geometry 
which motivate some of the tree theory are discussed on an elementary level.
Readers impatient to see technical definitions should skip ahead to
Section~\ref{ss de} and read those before continuing, then read Section~2.2
for a prototypical application and Section~2.3 for statments of all the
theorems.  Section~3 applies these to branching models and discusses
several of the ways that analysis of branching random walk may be applied
to models of disarate physical phenomena.  Section~4 applies the tree
theory to the geometry of random Cantor-like sets.
Finally, Section~5 mentions some problems of interest that are 
internal to the theory of tree-indexed processes.

\subsection{Branching models} \label{ss br}

The {\em simple} or {\em Galton-Watson} branching process models the
family tree of descendants from a single progenitor.  This individual has
a random number of children (possibly zero), each of which in turn has
a random number of children, and so on, with each of these random numbers
being independent picks from the same {\em offspring distribution}.  
The resulting random tree was studied in the previous century by 
Bienaym\'e, by Galton and Watson, and subsequently by others; see 
Heyde and Seneta (1977) for some of the history.

Many variants have been considered.  The {\em multitype} process 
separates individuals into different types (usually finitely many),
where each type, $1 , \ldots , m$, has a different distribution for 
the vector $(X_1 , \ldots , X_m)$ of the numbers of offspring it
will have of each of the $m$ types.  Instead of varying according
to the type of the parent, the offspring distribution may vary 
with each successive generation.  A {\em branching process in a
varying environment} (BPVE) has, instead of a single offspring
distribution,  a sequence of distributions, and all individuals in
generation $n$ have numbers of offspring that are independent
draws from the $n^{th}$ offspring distribution.  
The genetic applications of these models are obvious.  Interest from
another angle was sparked in the 1930's and 40's by the study
of cosmic ray cascades, electron-photon cascades, and of 
nuclear chain reactions.

Suppose now that each individual is born at a specified location,
displaced from its parent by a random vector, and that these
vectors are independent and identically distributed.  If $v$
is a vertex of the tree (i.e.\ an individual), let $X(v)$ denote
its displacement from its parent and let $S(v)$ be its location,
which is the sum of $X(w)$ over all ancestors of $v$ including
$v$ itself.  This process is called a {\em branching random walk}.
Branching random walks model many physical phenomena and the
study of their properties is far from exhausted.  Geneticists
and population biologists use branching random walks to model
dispersion of species, of genes and of infectious diseases.  An
example along these lines is discussed in detail in Section~2.2.  The 
remainder of this section is devoted to describing the various applications 
of branching random walks to other probability models and 
the mathematical questions that these generate.  

Interpreting the IID displacements as time lags gives a model called 
{\em first-passage percolation}.  The basic question is: what generation 
is reached by what time?  This was originally intended to model the
diffusion of liquid in a porous material (the graph being a Euclidean
lattice rather than a tree).  To model a chain reaction, one would
naturally use a tree whose vertices represented the events in the
chain reaction; one could also model the progress of a parallel
computation by first-passage percolation on the decision tree.
First-passage percolation may be applied to the characterization of
a random set known as {\em diffusion-limited aggregation}.  This is a
model for the growth of a cluster of particles in which each 
subsequent particle sticks to the existing cluster at a random location,
distributed according to the hitting measure of a random walk started
at infinity (Barlow, Pemantle and Perkins 1993). 
 
Reinterpreting the displacements as resistances of segments of wire 
gives a random electrical network that is mathematically equivalent to
a random walk on a tree in which the transition probabilities are
themselves random (a {\em random walk in a random environment});
see Doyle and Snell (1984) for the connection between random walks
and electrical networks.  The random walks in random environments
are in turn equivalent to certain {\em reinforced random walks}, in 
which the probability of a transition increases each time the transition
is made (Pemantle 1988).  Reinforced random walks are models for
learned behavior, and while trees are not the natural graphs on which
to run RRW's, they are to date the only graphs on which RRW's are
at all tractable (with the possible exception of some essentially
one-dimensional graphs).  

The IID displacements may be interpreted as energies.  This results in 
a thermodynamic ensemble having density $e^{-\beta H}$ with respect 
to product measure, where $H$ is the energy of a state.  Lyons (1989)
discusses an Ising model, in which a state is an assignment of $+1$ or
$-1$ to each vertex and $H$ is the sum of all edge energies.  Derrida 
and Spohn (1988) discuss a polymer model in which the states are 
paths of length $n$ in a regular tree of depth $n$ and $H$ is the
sum of energies along the path.  In either case,
an exponentially small probability (with repect to the reference measure) 
of an underaverage value of $H$ can greatly influence the partition 
function, and hence information about the behavior of a typical
element of the ensemble.  Since $H$ is determined from partial sums 
of IID random variables, one is led again to the extremal value theory 
of branching random walks.  In these models the tree structure is 
not completely natural, but is instead an approximation to the mean-field
limit in high dimensions; see Derrida and Spohn (1988) for a fair amount
of justification of the model.  Interpreting the IID displacements
as intensity factors of rainfall gives the cascade model for 
spatial distribution of rainfall studied by Gupta and Waymire (1993).
These random, stochastically self-similar, hierarchical spatial
distributions have been studied in other contexts by Kahane and 
Peyri\`ere (1976), by Waymire and Williams (1994), and others.
Finally, we will see in Section~4 how branching random walks may be
used to encode and solve problems in fractal geometry.

One basic question that arises in all these applications is the
extremal value question.  If the locations are one-dimensional, 
one might ask for the maximal displacement likely to occur 
in generation $n$ as a function of $n$.  To the first order, this 
is linear in $n$ and the method of computing the constant is
well known; this will be discussed at length in Section~3.  The deviation 
 from this was computed by Bramson (1978) and Derrida and Spohn (1988).
Another kind of extremal behavior is to ask whether there is an
infinite line of descent which exhibits a property which has probability
zero for any fixed line of descent.  For example, is there a line of
descent staying within a bounded region (Benjamini and Peres 1994b)?  
Is there a line going to infinity at a specified rate (Pemantle and Peres 
1994)?  The classical questions about branching processes (time to 
extinction, rate of growth) may also be phrased in terms of the
extremal value question, though a discussion of this would be 
too far afield.

The classical method for studying branching models is via generating 
functions.  Generating functions for the population at generation $n$
may be written exactly in terms of the generating function for
the offspring of each individual.  This method is powerful, but
often breaks down when events are weakly dependent rather than 
independent.  By contrast, the tree-indexed method proceeds as follows.
First, calculate the probabilities of seeing various things along a single
line of descent.  The probability of a single line staying in
a given region or escaping to infinity at a given rate is a classical
computation since the increments are IID.  Multiplying the expected
size of generation $n$ by the probability of a given behavior gives
the mean incidence of that behavior.  The probability of observing
the behavior is bounded by the mean incidence, but may be less; a 
second moment computation will distinguish between these cases.  In
other words, the mean incidence tells you the one-dimensional distributions
of a random field indexed by the boundary of the tree, and the
second moments give you enough information about the joint distributions
to get probability bounds.

Results on branching models are worked out in Section~3, with a
prototypical argument previewed in Section~2.2.  While the sharpest 
results on tree-indexed processes are all stated in terms of 
potential theory, the previous paragraph should serve as a guide
to the structure of the arguments for non-experts in potential theory.

\subsection{Random sets with stochastic self-similarity} \label{ss sss}

In Section~4, it will be shown how to make a correspondence between
paths in an infinite homogeneous tree and points in Euclidean space,
which preserves the potential-theoretic structure.  Consequently,
questions about random subsets of Euclidean space may be analyzed
in terms of the corresponding random trees.  In particular, certain
stochastically self-similar sets correspond to well understood 
random trees, such as Galton-Watson trees, making knowledge
especially easy to transfer.  Self-similar and stochastically self-similar
sets are usually {\em fractals}, meaning that the have a non-integral
dimension.  These have been widely studied in the last 20 years, both
as complex mathematical objects and as visually beautiful objects whose
scale-invariance captures some intriguing aspects of natural law; see
Falconer (1985) for a mathematical introduction and consult the science
section of your local bookstore for pretty pictures.  

Consider the following Cantor-like set.  Let $A_1 , \ldots , A_N$
be a collection of subcubes of the $d$-dimensional unit cube.  
We allow $N$ and $A_1 , \ldots , A_N$ to be random but require
that their law $\mu$ concentrate on collections with disjoint interiors.  
Let $C_1 = \bigcup_k A_k$, or in other words, throw out everything 
not in one of the sets $A_j$.  Apply this recursively
to each $A_j$: choose a collection of subcubes $A_{j,1} , \ldots ,
A_{j,N_j}$ independently from the image of $\mu$ under the similarity 
that maps $[0,1]^d$ to $A_j$, and throw out everything in 
$A_j \setminus \bigcup_r A_{j,r}$.  The limiting set $C$ is stochastically
self-similar in an obvious sense.  Familiar examples are as follows.
If $d=1$ and $\mu$ is a point mass at the collection $\{ [0,1/3] ,
[2/3 , 1] \}$ one gets the usual (deterministic) Cantor set.  If
$d=1$ and $\mu$ picks $\{ [0,a] , [b,1] \}$ with $(a,b-a, 1-b)$ having
Dirichlet ($1/2,1/2,1/2$) distribution, then $C$ is distributed as the zero
set of a Brownian bridge\footnote{Many other distributions for $(a,b)$ 
generate the Brownian zero set as well; the present example may be
found in Perman, Pitman and Yor (1992).}.
If $d=2$ and $\mu$ gives probability
$p^k (1-p)^{9-k}$ to every subcollection of size $k$ of the 
partition into 9 squares of side $1/3$, then one gets the 
so-called canonical curdling process studied by Chayes, Chayes and 
Durrett (1988) and Dekking and Meester (1990).  

Hawkes (1981) computes dimensions of a large class of such sets.
Graf, Mauldin and Williams (1988) compute precise Hausdorff gauges.
Chayes, Chayes and Durrett (1988) and Meester (1990) 
discuss connectivity properties, but these problems seem
to be hard and no general criteria are known.  The approach
carried out in Section~4 is to determine properties such as dimension
by establishing close connections between the random sets and 
the representing trees, then to use known facts or relatively easy 
theory to analyze the trees.  An advantage to this method is 
that the dimension may be bounded below without exhibiting a
measure meeting the appropriate regularity condition.  In some 
sense, the methods used by Hawkes, by Graf-Mauldin-Williams, 
by McMullen (1984), and earlier by Carleson and Frostman 
are all based on the idea of a tree representation.  

Another question about a random set is its intersection properties.  
For example, Can you tell when two random sets have positive 
probability of intersecting? Two independent Brownian motions 
intersect in dimensions less than 4, while three or more intersect 
only in dimension 2.  A complete characterization exists of 
which sets intersect Brownian motion with positive probability 
(Kakutani 1944), but only recently was it determined which sets 
have a common intersection with two Brownian motions.  
Fitzsimmons and Salisbury (1989) solve this problem using 
classical potential theory, settling a conjecture 
of Evans and of Tongring, while Peres (1994a)
has a much simpler proof translating the problem to trees.  This approach
also shows how to compute the drop in the dimension of a set when intersected
with various stochastically self-similar sets including the range of a
Brownian motion.  Peres (1994b) describes several other applications 
resulting from translating geometric questions about 
Brownian motions to trees.  Related to these results are two
theorems of Marstrand, showing that positive one-dimensional capacity 
is sufficient for a set to intersect a random line with positive 
probability (the converse fails but not by much) and that
the dimension of the intersection is, generically, one less than 
the dimension of the original set.  A derivation of the latter
 from the former may be established by tree methods.  
Section~4 discusses these results in more detail.  

\subsection{Other motivations} \label{ss rw}

A significant part of the motivation for studying the potential theory of trees
came from random walks on trees.  The geometry of a Riemannian manifold
can be analyzed in terms of the behavior of Brownian motion on the 
manifold (see Ledrappier 1988).  Negatively curved manifolds may be
discretized so that Brownian motion on the manifold corresponds to a 
random walk on an embedded tree.  Symmetric spaces give rise to
periodic embedded trees, manifolds of negative curvature bounded away from zero
have embedded trees of exponential growth, and so on.  Lyons (1993)
discusses behavior at infinity of random walk on periodic trees, while
Lyons, Pemantle and Peres (1995) discuss the randomized counterpart,
where the tree is Galton-Watson.  Conditions for the recurrence or 
transience of random walks on arbitrary trees (in terms of capacities)
were obtained by Lyons (1990) and by Benjamini and Peres (1992b).

Homogeneous trees are Cayley graphs of free groups, and random walks on
trees qua Cayley graphs have been studied by many people; 
see Mohar and Woess (1989) for
some references to studies of random walks on homogeneous trees.
Sawyer (1978) proposes a random walk on a tree as a model for the
dispersion of genetic types along a river system.  The model is 
very rough, but the spectral and boundary theory there is shown
to answer natural questions about the distribution of types. 
While random walks on trees do not 
constitute tree-indexed processes (for which the tree should be the index
set, not the range space), they share the same techniques.
For instance, the classification of recurrence/transience of random walks
on a tree in Benjamini and Peres (1992b) results in the same capacity
criterion as for a certain set to be polar in Pemantle and Peres (1995a);
this is not entirely a coincidence, and a more explicit connection 
is made at the end of the final section of this article.  At any rate,
random walks on trees and tree-indexed random walks have cross-fertilized 
each other enough to warrant mention here of the latter.

In addition, the generalization of branching random walks in which the
branching part is deterministic and given by an arbitrary tree has
been studied for its own interest.  The first mention of this is by
Joffe and Moncayo (1973), although it was not wholeheartedly pursued 
until Benjamini and Peres (1994a), having been generalized meanwhile 
to Markov chains indexed by trees in Benjamini and Peres (1992a and 1994b). 

\section{Technical overview} 

\subsection{Definitions} \label{ss de}

\noindent{A} tree is a connected, undirected graph with no cycles.
All trees are assumed as well to be locally finite (i.e.\ finitely
many edges incident to each vertex) and to have a distinguished
vertex known as the root.  The name used most often for a generic
tree is $\Gamma$ and its root will most often be denoted $\rho$.  
The name $\TB$ is reserved for the infinite $b$-ary tree, in which
each vertex has $b$ children (neighbors at greater distance from
the root).  The notation $x \in \Gamma$ will be used for
``$x$ is a vertex of $\Gamma$'' since no confusion results.  
Let $|x|$ denote the number of edges in the path connecting
$x$ to the root and let $\Gamma_n$ denote the set $\{ x \in \Gamma :
|x| = n \}$ of vertices in the $n^{th}$ level or generation
of $\Gamma$.  For vertices $x,y \in \Gamma$, define $x \leq y$
to be the relation that holds if $x$ is on the path from
$\rho$ to $y$, and let $x \wedge y$ denote the greatest lower
bound of $x$ and $y$ (i.e.\ the vertex at which the paths from
$\rho$ to $x$ and $y$ diverge).  For reasons to be seen 
shortly, trees are usually assumed to have uniform height.
A tree of height $N < \infty$ has uniform height if all
its leaves (vertices without children) are at level $N$; a tree
of uniform height $\infty$ has no leaves at all.  The boundary
$\partial \Gamma$ of a tree $\Gamma$ of height $N \leq \infty$ is 
the set of self-avoiding paths of length $N$ starting from the root.
If $\Gamma$ has height $N$ but not uniformly, then $\partial \Gamma$ 
contains paths through only those vertices with descendants at 
level $N$; since $\partial \Gamma$ is of fundamental interest,
vertices with no descendants at level $N$ become irrelevant,
whence the assumption that there aren't any.  Extend the symbol
``$\wedge$'' to $\partial \Gamma \times \partial \Gamma$ by
letting $x \wedge y$ denote the greatest vertex of $\Gamma$ 
contained in both $x$ and $y$. 

Sometimes the trees are random, the most common type of
random tree being a Galton-Watson tree.  This is the
family tree of a branching process in which each individual
has a random number of children and all these numbers
are IID.  The usual notation for Galton-Watson trees 
is in effect: $f (z) = \sum a_n z^n$ is the offspring generating
function, where $a_k$ is the probability of having $k$
children and $f' (1) = \sum k a_k$ is the mean number of
offspring per individual.  The law of a Galton-Watson
tree with offspring generating function $f$ is denoted
$GW^f$, or just $GW$ when $f$ is clear from context.

Let $\mu$ be a probability distribution on a measure space $S$,
often taken to be the real numbers, and let $\{ X(v) : v \in \Gamma \}$
be a collection of IID random variables, indexed by the vertices
of $\Gamma$, having common law $\mu$, and defined on the measure
space $(\Omega , \F , \P)$.  Give $S$ the discrete topology, 
in which all sets are open (though not necessarily measurable),
and give $S^N$ the product topology, which is discrete unless
$N = \infty$.  Let $B \subseteq S^N$ be any measurable closed set.  
Define an event $A$ depending on $\Gamma$, $B$ and $\Omega$ by
$$A(\Gamma ; B ; \Omega) = \left \{ \exists (\rho , v_1 , v_2 , \ldots)
   \in \partial \Gamma \,:\, (X(v_1) , X(v_2) , \ldots) \in B \right \} .$$
In other words, $A$ is the event that there is some path for which
the sequence of values of the $X$'s lies in the prescribed set, $B$.
The quantity $\P (A(\Gamma ; B ; \Omega))$, which depends on
$\Omega$ only through $\mu$ is denoted $P(\Gamma ; B ; \mu)$, or
when $\mu$ is understood, just $P(\Gamma ; B)$.  Viewing the
probability space $\Omega$ as defining a $S^N$-valued random field
on $\partial \Gamma$, the first natural question is which sets
$B$ are ``hit'' by the random field (intersect its range with
positive probability).  Sets $B$ for which $P(\Gamma ; B) = 0$ 
are called {\em polar} by Evans (1992) and thus the classification of
sets as polar or nonpolar becomes the primary object of study.

A special case is when $B$ is the product set 
$$[0,a_1] \times [0,a_2] \times \cdots .$$
This is called {\em Bernoulli percolation} by Lyons (1992).
In this case one may imagine killing vertices randomly
and independently, killing a vertex in generation $n$ with
probability $1 - a_n$; then $B$ is the set of paths all of whose
vertices remain alive.  The independence makes this case
easier to analyze, and the first and sharpest theorems
were obtained here. 

\subsection{The basic idea: second moments} \label{ss 2mm}

In order to illustrate the use of potential theory, I devote this section
to working out a simple branching random walk example.  Consider a
flower germination model, beginning with a single individual, which
sends out during the course of its lifetime $b$ spores, the locations
of which are displaced from the parent by vectors that are IID$\,\,\sim\mu$.  
I have assumed for simplicity that the branching is deterministic.
In the notation of the preceding section, $\Gamma = \TB$, $S = \R^2$
and $N = \infty$.  Fix a region $G \subseteq \R^2$ representing
hospitable terrain, and suppose that spores alighting outside of $G$
fail to germinate.  Such models and variants thereof can be found in 
Levin et al (1984) and Bergelson et al (1993), among other places.

Let $B = \{ (x_1 , x_2 , \ldots ) :  \sum_{i=1}^n x_i \in G 
\mbox{ for every } n \}$.  Then $A (\Gamma ; B)$ is the event of
nonextinction of this flower's family tree.  Suppose $G$ is
a nice set: a compact closure of a connected domain.  Let $\pi_n (B)$
denote the projection of $B$ onto the first $n$ coordinates, i.e., those
paths staying in $G$ for at least the first $n$ steps.  The probability
$P (\Gamma ; B)$ is the decreasing limit of probabilities 
$P (\Gamma ; \pi_n (B))$.
Now for a single line of descent, the probability $\mu^n (\pi_n (B))$ of
staying inside $G$ for $n$ steps is easy to estimate: it is asymptotically
a constant multiple of $\lambda^n$ where $\lambda$ is the maximal 
eigenvalue of the region $G$.  Let $W_n$ be the number of survivors in 
generation $n$, so $\E W_n = (b \lambda)^n$ up to a constant factor.  
Obviously the process must die out when $\lambda < 1/b$.  To show that the
process may survive when $\lambda > 1/b$, we show that 
$\E W_n^2 / (\E W_n)^2$ is bounded by some constant $C$ independent 
of $n$.  This directly implies that $\P (W_n > 0)^{-1} \leq C$ for
all $n$, and hence $P(\Gamma ; B) \geq C^{-1}$
(see Aldous 1989 for some other uses of this implication).  

To compute $\E W_n^2$, sum over pairs $(v,w) \in \Gamma_n^2$
the probability that both $v$ and $w$ have lines of
ancestry staying completely inside $G$.  Clearly this probability
depends only on $n$ and $|v \wedge w|$.  In fact, conditioning
on $\sum_{z \leq v \wedge w} X(z)$ shows it to be bounded
above by a constant multiple of $\lambda^{2n - |v \wedge w|}$.
Thus we may write 
$${\E W_n^2 \over (\E W_n)^2} \leq c \sum_{v,w \in \Gamma_n}
   \lambda^{2n - |v \wedge w|} (b \lambda)^{-2n} = c \int\int
   \lambda^{- |v \wedge w|} \, d\mu (v) \, d\mu (w) $$
where $\mu$ is the uniform measure on $\Gamma_n$.  An easy computation 
shows this is finite when $\lambda > 1/b$, which finishes the demonstration.

As a preview, consider what would have happened if $\Gamma$ were
not a homogeneous tree.  Let $K(v,w)$ denote $\lambda^{- |v \wedge w|}$.
If $K$ has a finite integral against the product uniform measure,
then the same argument shows that a branching random walk
indexed by $\Gamma$ has a line staying in $G$ with positive
probability.  In fact, making $W$ a weighted sum of indicator
functions of lines of descent staying in $G$ shows that the
measure one integrates against need not be uniform.  Furthermore,
any measure $\mu$ on $\partial \Gamma$ with 
$$\int \int K(x,y) \, d\mu^2 < \infty$$
projects to a measure $\mu_n$ on each $\Gamma_n$ for which
the integrals of $K$ against $\mu_n^2$ are bounded.  Thus one
obtains the result: if $\partial \Gamma$ supports a 
probability measure $\mu$ with $\int \int K(x,y) \, d\mu^2 < \infty$,
then $P (\Gamma ; B) > 0$ .  Restating the hypothesis of this
result in the language of potential theory gives the Basic Theorem of
the next section.

\subsection{Statements of theorems} \label{ss th}

Since the notions of energy and capacity are fundamental to the results
surveyed here, I include a brief discussion.  A few 
definitions and examples are no substitute for familiarity, 
so the reader is referred to Carleson (1967, Chapters I - IV), 
or to Falconer (1985, Chapter 6) for geometric facts 
about metric capacity.  The relation between capacity and dimension
is that the capacity of a set $A$ in gauge $x^{-\alpha}$ will
be positive if $dim (A) > \alpha$ and zero if $dim (A) < \alpha$.

Given a probability measure $\mu$ on a metric space and given
a monotone function $g$ on the positive reals tending to infinity 
at zero, the energy of $\mu$ with respect to $g$ is defined by
$$\cE_g (\mu) = \int \int g(d(x,y)) \, d\mu (x) \, d\mu (y) ,$$
where $d(x,y)$ denotes distance.  The $g$-capacity of a set $A$ 
is defined by
$$ \Cap_g (A) = \left [ \inf \{ \cE_g (\mu) : \mu (A) = 1 \} \right ]^{-1} .$$ 
There is a natural class of metrics on $\partial \Gamma$ gotten
by letting $d(x,y)$ be any function of $|x \wedge y|$ that
decreases to zero as $|x \wedge y| \rightarrow \infty$.  The
notion of metric energy and capacity on $\partial \Gamma$ for 
these metrics may be formulated directly in terms of functions 
$f : \Z^+ \rightarrow \R^+$ that increase to infinity:
$$\cEB_f (\mu) = \int \int f(|x \wedge y|) \, d\mu (x) \, d\mu (y)$$
for probability measures $\mu$ on $\partial \Gamma$, while
$\Capb_f (A)$ is the reciprocal of $\inf \{ \cEB_f (\mu) : \mu (A) = 1
\}$, as before.  

In this language we may restate the result from the previous section,
stated and proved in Pemantle and Peres (1995a) but already implied by
Lyons (1992).   

\vspace{1ex}\noindent{\bf Basic Second Moment Theorem:~~}{\it
Let $\Gamma, N, S, B, \mu$ and $\{ X(v) : v \in \Gamma \}$ be as
in Section~\ref{ss de}.  Let $W$ denote the set of vertices $v$
such that $(X(v_1) , X(v_2) , \ldots , X(v)) \in \pi_{|v|} (B)$,
where $\rho, v_1 , v_2, \ldots , v$ is the path connecting the root
to $v$ and $\pi_k$ is the projection onto the first $k$ coordinates.
Suppose there is a positive, nondecreasing function 
$f: \Z^+ \rightarrow \R^+$ such that for any two vertices 
$v,w \in \Gamma$ with $|v \wedge w| = k$, 
\begin{equation} \label{eq QB}
\P (v,w \in W) \leq f(k) \P (v \in W) \P (w \in W)
\end{equation}
Then 
$$P (\Gamma ; B) \geq \Capb_f (\partial \Gamma) .$$
}\\[1ex]

\noindent{\em Remark}: Usually, when a second moment 
(also known as $L^2$) method is used, 
there is a question as to whether the result is sharp.  
If you followed the argument in the previous section, you will
notice that the property of $W$ stated in~(\ref{eq QB}) is 
enough to imply the conclusion regardless of whether any
variables $\{ X(v) \}$ underlie the definition of $W$.
One cannot expect sharpness without using further properties
of $W$, which will now be explored.

Typically, the second moment method shows some property to hold 
if a set satisfies $\Capb_f (A) > 0$, while a simpler first moment 
estimate shows the converse to hold if $A$ has 
zero Hausdorff measure in gauge $f$; this leaves
a small gap\footnote{In all nontrivial cases for which I know
the resolution of the gap, the capacity criterion is
sharp, not the measure criterion.  See Kahane (1985) for
some instances of the gap, e.g.\ Theorem~5 on page 246 and~(5)
on page~236; see Shepp (1972) for a resolution of the gap in one
case, in favor of the capacity criterion.}.  One circumstance in which
the Basic Theorem is sharp, up to a factor of 2, is when 
$B$ is Bernoulli.  Recall that $B$ is Bernoulli if 
it is a product set, $B = \prod_n [ 0,a_n ]$.
In this case $f(k)$ may be taken to be $\P (v \in W)^{-1} =
\mu^k (\pi_k (B))$ for $v \in \Gamma_k$.  This is clearly
the least $f$ can be (take $w = v$ in~(\ref{eq QB})).  

\vspace{1ex}\noindent{\bf Sharp Bernoulli Theorem (Lyons 1992):~~}{\it
If $B$ is Bernoulli and $f(k) = \P (v \in W)^{-1}$ for any $v \in 
\Gamma_k$, then 
\begin{equation} \label{eq sharp QB}
2 \Capb_f (\partial \Gamma) \geq P(\Gamma ; B) \geq \Capb_f
   (\partial \Gamma) .
\end{equation}
}\\[1ex]

Cases where $f(k)$ may be taken as $C \P (v \in W)^{-1}$ are called
{\em quasi-Bernoulli}.  Here too, the gauge function is as 
small as possible (constant multiples being inconsequential)
and a converse is conjectured; see Section~4.  
   
Many potential theoretic results from Markov process theory are of
the form: A Markov process hits a set $A$ with positive probability
if $A$ has positive capacity in a certain gauge (determined from the
Green's function of the process).  The most famous of these is due to
Kakutani (1944) and is sharp: Brownian motion in $\R^d$, $d \geq 3$
hits a set $A$ with positive probability if and only if $\Cap_g (A) > 0$,
where $g(x,y) = |x-y|^{2-d}$.  The previous results were dual to
this, in that they gave capacity conditions on $\Gamma$ rather than
on $B$.  Here is a direct tree-indexed analogue of Kakutani's theorem.

\vspace{1ex}\noindent{\bf Dual Second Moment Theorem:~~}{\it
Suppose $\mu$ is the uniform distribution on the
set $\{ 1 , 2 , \ldots , b \}$.  Let $\Gamma = \TM$ be the 
homogeneous $m$-ary tree, that is, a tree where the root has
$m$ children and each other vertex has $m+1$ neighbors, those
being the parent and $m$ children.  Observe that the closed
set $B$ is naturally encoded as a subset $\overline{B}$
of the boundary of the $b$-ary tree.  In this notation, 
$$2 \Capb_f (\partial \overline{B}) \geq P (\Gamma ; B) 
  \geq  \Capb_f (\partial \overline{B}) , $$
where $f(k) = k$ if $m = b$ and $f(k) = (b/m)^k$ if $m < b$.
When $m > b$, $P(\Gamma ; B) > 0$ for all nonempty $B$.
}\\[1ex]

Note that the Dual Second Moment Theorem is sharp (``if and
only if'') but at the expense of restricting to homogeneous
trees, which is analogous to restricting to Bernoulli sets.
This theorem was first proved by Evans (1992) with a factor
of 16 instead of 2 and by Lyons (1992) with a factor of 4.
The proof with a factor of 2 follows from the methods 
of Benjamini, Pemantle and Peres (1993).  

The next two theorems give conditions for one tree to have all
the polar sets that another tree has.  Such comparisons are
useful because only in the case of homogeneous trees are the
polar sets easy to compute.  Say that a tree $\Gamma^1$ is
at least as polar as $\Gamma^2$ if every polar set for $\Gamma^2$
is a polar set for $\Gamma^1$, and call two trees ${\em
equipolar}$ if they have the same polar sets. 
A tree $\Gamma$ is called {\em spherically symmetric} if each 
vertex $v$ in $\Gamma_{n-1}$ has precisely $f(n)$ children
for some function $f$.

\vspace{1ex}\noindent{\bf Comparison Theorem:~~}{\it
Suppose $\Gamma$ is spherically symmetric 
and let $\Gamma'$ be any tree with $|\Gamma_n'| \leq |\Gamma_n|$
for all $n$.  Then $P (\Gamma' ; B) \leq P (\Gamma ; B)$
for any $B$ and $\mu$.  
}\\[1ex]

As mentioned earlier, the tree $\Gamma$ may itself be random.
The following is a ``universality
class'' theorem for Galton-Watson trees, saying that in the
finite variance case, trees with the same mean growth
are equipolar and thus essentially the same from a tree-indexed
process point of view.

\vspace{1ex}\noindent{\bf Equipolarity Theorem:~~}{\it
Let $GW_1$ and $GW_2$ be the Galton-Watson measures, 
corresponding to two offspring distributions with 
the same mean $m > 1$ and each having finite variance.
(Zero variance is allowed in the
case that $m$ is an integer.)  Then for $G_1 \times G_2$-almost
every $(\Gamma^1 , \Gamma^2)$ there exist almost surely constants 
$0 < c < C < \infty$ depending on $\Gamma^1$ and $\Gamma^2$
such that for all $\mu$ and all sets $B$, 
$$c P (\Gamma^1; B) \leq P(\Gamma^2 ; B) \leq C \P (\Gamma^1 ; B) .$$
In particular, $\Gamma^1$ and $\Gamma^2$ are equipolar.
If instead $\Gamma^1$ has infinite variance, then
the above does not hold, and in fact $\Gamma^1$
has strictly more polar sets.
}\\[1ex]

The fact that the trees are not equipolar when one offspring
variance is infinite should provide some resistance against
the notion that the Equipolarity theorem is obvious.  For
more evidence, consult Graf et al (1988), 
wherein it is shown that Galton-Watson trees do not behave 
the same as deterministic trees of the same mean with
respect to Hausdorff measure.  The Comparison Theorem is
 from Pemantle and Peres (1994) and the Equipolarity Theorem
is from Pemantle and Peres (1995b) and Pemantle (1993).

\section{Applications to branching models}

For any tree-indexed process $(\Gamma , B , \Omega)$ whose
state space is a group, one may define
partial sums (or products in the non-abelian case) by
$$S(v) = \sum_{\rho < w \leq v} X(v)$$
for each $v \in \Gamma$.  If the state space is $\R$, 
define the extremal values by
$$M_n = \max \{ S(v) : |v| = n \} .$$
Questions about $M_n$ have arisen in the contexts of
random distribution functions (Dubins and Freedman 1967), directed
polymers and partition functions for high-dimensional limits 
of random fields (Derrida and Spohn 1988), 
branching random walks (Bramson 1978), a more general 
``Markov branching random walk'' (Karpelevich et al 1993), 
reinforced random walks (Pemantle 1988), random walks in 
random environments (Lyons and Pemantle 1992), as well as 
indirectly in the study of explosions in first-passage percolation 
(Pemantle and Peres 1994).  

Consider a reasonably simple case.  Suppose that $\Gamma = \TB$
is a homogeneous tree and the common distribution $\mu$ of the
$\{ X(v) \}$ is bounded.  If $\alpha > m = \E X( \rho )$ then
a standard large deviation estimate yields
\begin{equation} \label{eq alpha}
\P (S(v) \geq \alpha |v|) = (C_\alpha + o(1)) |v|^{-1/2} 
   u(\alpha)^{|v|} ,
\end{equation}
where $u (\alpha )$ is the rate function.
Clearly, if $u (\alpha) \leq b^{-1}$ then $M_n \leq \alpha n$ 
with high probability.  Pick $\alpha_0 \in (0,\infty)$ to be
the infimum of $\alpha$ for which $u(\alpha) \leq b^{-1}$.  
Is $\alpha_0$ the correct limit of $M_n / n$?
To complete the picture, one must show that $M_n \geq 
(\alpha_0 - \ee) n$ for any $\ee > 0$ and sufficiently large $n$.  
The first proof is due to Hammersley (1974), who proved
convergence in probability of $M_n / n$; here is a tree-indexed
proof of almost sure convergence.

Let $B$ be the set 
$$\{ (x_1 , x_2 , \ldots ) : \sum_{i = (j-1)k+1}^{jk} x_i \geq
   (\alpha_0 - \ee) k \mbox{ for all } j \} .$$
It is easy to verify quasi-Bernoullicity, hence $\P (\Gamma ; B)
\geq \Capb_f (\partial \Gamma)$ where $f (|v|) = C \P (v \in T)^{-1}$.
The choice of $\alpha_0$ guarantees that for fixed $\ee$ and 
large enough $k$, $f(k) \leq (2-\dd)^k$ for some $\dd > 0$.
The binary tree has positive $(2-\dd)^k$-capacity for every $\dd > 0$, 
so $P (2^\infty ; B) > 0$.  On this event, $\liminf M_n / n \geq
\alpha_0 - \ee$.  But the event $\liminf M_n / n \geq \alpha_0
- \ee$ is a tail event in the $\{ X(v) \}$, so it has probability\
one, and since $\ee > 0$ is arbitrary, this yields $\lim M_n / n = \alpha_0$.  

As we discuss this argument, let us compare the result in varying
degrees of generality to a string of such results proved in the
last 25 years.  Dubins and Freedman (1967) consider the case
where $\Gamma$ is a binary tree and the $\{ X(v) \}$ are
Bernoulli ($p$).  They observe that when $p > 1/2$, there is
a path with only finitely many 0's, and ask what the maximum
density of 1's along a path is for $p < 1/2$.  They solve
the easy half, using Borel-Cantelli to show that the density 
can be at most the value $\alpha_0$ above, leaving the other
direction as an open problem (page 207).  From the modern perspective, 
it is striking that this question was difficult to settle at the time!

The earliest solution I know of in the mathematics literature was by 
Kingman (1975), in the context of first-birth times.  He allows an 
arbitrary common distribution on the positive reals for $\mu$ and 
considers the tree $\Gamma$ to be the family tree of a Galton-Watson
branching process with mean growth $m$, conditioned on survival.  
Kingman's proof is quite specific to this particular 
problem, relying on exact computation of 
$$\E \left [ \sum_{|v| = n} S(v) e^{- \theta S(v)} \right ] \; .$$
The result is that $M_n / n \rightarrow \alpha_0 (m)$, where 
$\alpha_0 (m)$ is gotten by replacing the 2 in~(\ref{eq alpha}) by $m$.
Biggins (1977) allows $X(v)$ to take negative values, provided 
the moment generating function exists in a neighborhood of 0, and 
provides a simpler proof based on finding a supercritical branching
process in a tree derived from $\Gamma$ by looking only at vertices
in levels $0, k, 2k, \ldots$; essentially 
the same proof appears in Pemantle (1988).  Biggins also
computes the asymptotic number of $v$ in level $n$ for which
$S(v) \geq \alpha n$ when $\alpha < \alpha_0$.  Lyons and
Pemantle (1992) are the first to provide a proof via 
quasi-Bernoullicity; see also Kesten (1978).  
For the particular result $M_n / n \rightarrow
\alpha_0$, Biggins' proof is the simplest and best.  The main
advantage in the potential theoretic proof is its wider scope, allowing
for completely general trees or general growth rates of the $M_n$.
The result immediately generalizes from Galton-Watson trees of mean
growth $m$ to any tree that has positive capacity for gauges
$f(n) = r^n$ when $r < m$ 
and not when $r > m$; these trees include {\em periodic} trees of
mean growth $m$ (defined, for example, in Lyons 1993) as well as 
most reasonably small perturbations of homogeneous and 
Galton-Watson trees.   

If $S(v)$ is interpreted as (the negative of) an energy function, 
as in Derrida and Spohn (1988),
then one is interested in the behavior of the partition function
$Z_n \deq \sum_{|v| = n} \exp (S(v))$, and there should be
a limit $(\log Z_n) / n \rightarrow \beta$.  The contribution to
$Z_n$ from vertices $v \in \Gamma_n$ with $S(v) \approx \lambda n$
is $\exp [(\lambda - \phi (\lambda)) n ]$, where $\phi (\lambda)$ is the
large deviation rate for the average variables that are 
IID$\sim\mu$ to exceed $\lambda$.  In other words, the rate
is the same as if the values of $S(v)$ were independent for
$v \in \Gamma_n$\footnote{Lyons and Pemantle (1992) give a short proof,
but in fact this was proved by most of the authors cited above:
Kingman, Biggins, Derrida and Spohn, and by Kahane and 
Katznelson (1990) and Waymire and Williams (1994) in the 
context of cascade spectra.}.  Taking 
\begin{equation} \label{eq dim}
\beta = \sup \{ \lambda - \phi (\lambda) : \lambda \leq \alpha_0 \}
\end{equation}
gives the correct limit: $(\log Z_n) / n \rightarrow \beta$ almost surely.
Observe for application in Section~4 that if $\exp (S(x \wedge y))$
defines a (random) metric on $\partial \Gamma$ and $\beta < 1$
then this shows the Hausdorff dimension of $\Gamma$ to be at most 
$\log b / |\log \beta|$.  In fact Theorem~4 of Lyons and Pemantle (1992) 
shows that this is exactly the dimension.

The Equipolarity Theorem gives universality results for the extremal
problem.  In the special case where $\Gamma$ is a binary tree and
$\mu$ gives measure one half to 0 and one half to 1, 
Bramson (1978) shows that the median value of $M_n$ differs from the
linear estimate by $K \log \log n$, where $K = (\log 2)^{-1}$.
The distribution is tight around its median as $n \rightarrow \infty$.  
Bramson, and later Derrida and Spohn (1988), carry out their analyses
on a binary tree, where an exact recursion reduces the problem
to a nontrivial analysis of the KPP equation.  Their results may
immediately be extended to the case of a general, mean-two, finite-variance
branching mechanism as follows.  Let 
$$B = B_{n,l} = \{ (x_1 , x_2 , \ldots) : x_n \geq l \} .$$
Apply the equipolarity theorem to Galton-Watson trees with 
mean two and finite variance (including the deterministic binary
tree) to conclude that there is almost surely some $c$ for which
the $\alpha$ quantile of $M_n$ on the Galton-Watson tree
is bounded above and below by the $1-(1-\alpha)/c$ and $\alpha / c$
quantiles respectively for $M_n$ on the binary tree $2^\infty$.  
This gives the new result that $M_n - K \log \log n$ is tight for 
Galton-Watson trees as a result of being tight for $2^\infty$.  
Similarly, the equipolarity theorem shows that survival with positive
probability in the flower germination model of Section~\ref{ss 2mm}
depends only on the region and the mean offspring, but not on 
the particular offspring distribution as long as it has finite variance.

The Comparison theorem turns out to be useful in the analysis of
first-passage percolation.  As mentioned before, first-passage percolation
has been used to describe a randomly growing subtree (Knuth (1973),
Aldous and Shields (1988), Barlow, Pemantle and Perkins (1994)).  
Andjel (personal communication) raises the question of how quickly
$\Gamma$ may grow and still have the minimum $m_n ( \Gamma ) \deq
\min \{ S(v) : |v| = n \}$ tend to 
infinity.  (This arises in a construction of an infinite particle system 
via a mapping from a collection of IID exponential random variables.)
The answer, given in Pemantle and Peres (1994), is that when $\Gamma$ is
symmetric with {\em growth function} $f$ in the sense that every vertex
in $\Gamma_n$ has precisely $f(n)$ children, then for increasing growth
functions and exponential random variables, $X(v)$,
\begin{equation} \label{eq explosion}
m_n \rightarrow \infty \mbox{ a.s.  if and only if } \sum f(n)^{-1} = \infty .
\end{equation}

The following sketch shows how the Comparison Theorem is instrumental
in obtaining a similar result in the case where $f$ is not
necessarily increasing.  For any $g$, let $\Gamma^g$ denote a 
spherically symmetric tree with growth function $g$.  Now choose 
a particular $g$, namely the pointwise greatest increasing integer
function for which $\prod_{k=1}^n g(k) \leq \prod_{k=1}^n f(k)$
for all $n$ (it is an easy exercise to verify the existence of such 
a $g$ and give other descriptions of it).  I will show 
that~(\ref{eq explosion}) holds with the $f$ on the right
replaced by $g$.  By definition, $|\Gamma^g_n| \leq |\Gamma_n|$ for 
all $n$.  The Comparison Theorem applied to the sets 
$B_{n,l} = \{ (x_1 , x_2 , \ldots ) : x_n \leq l \}$ 
implies that $m_n (\Gamma^g)$ is stochastically greater than 
$m_n (\Gamma)$.  This, along with~(\ref{eq explosion}) for
the increasing function $g$, proves that $\sup m_n < \infty$
whenever $\sum g(n)^{-1} < \infty$.  To prove the other
half, $\sum g(n)^{-1} = \infty \; \Rightarrow \; m_n \rightarrow 
\infty$, only a slight modification of
the proof for increasing functions is needed.

The Comparison Theorem is also used in Pemantle and Peres (1995b) to
prove one direction of the Equipolarity Theorem: almost every
Galton-Watson tree of mean growth $m$ has at least the same polar sets
as a certain spherically symmetric tree with more than $C_1 m^n$ 
children in generation $n$, showing the half of the Equipolarity 
Theorem that does not rely on finite offspring variance. 

\section{Dimensions and intersections of random sets}

The following correspondence between trees and Euclidean
space is vital to all of the applications in this section.  For convenience, 
I consider the unit cube $[0,1]^d$ rather than all of $\R^d$.  

Let $\C$ be the collection of closed binary subcubes 
of the unit cube, that is, cubes of the form 
$$\prod_{i=1}^d [j_i 2^{-n} , (j_i + 1) 2^{-n}]$$
where $n \geq 1$ and $2^n > j_1 , \ldots , j_d \geq 0$. 
The elements of $\C$ may be viewed as the vertices
of homogeneous $b$-ary tree ($b = 2^d$) in an obvious
way: a cube is a descendant of another if it is a subset.
Thus the root is the unit cube, and each cube has
$2^d$ children (descendants with no intervening lineage).
The identification of $\C$ with the vertices of $\TB$
induces a map $\phi$ from $\partial \TB$ to $[0,1]^d$, namely
$\phi$ of a sequence of cubes is the unique point in their
decreasing intersection.  Putting the metric 
$$ \mbox{dist}(x , y) := \sqrt{d} 2^{-|x \wedge y|}$$
on $\partial \TB$, it is clear that $\phi$ is continuous
and is in fact a contraction.  If $A$ is any closed subset 
of the unit cube, $\phi^{-1} [A]$ is a closed subset
of $\partial \TB$ and is therefore also the boundary of
a subtree of $\TB$.  Since $\phi$ is a contraction, it is
immediate that for any $g$, the metric capacities satisfy
$$\Cap_g (\phi^{-1} [A]) \geq \Cap_g (A) $$
(the pullback of any measure has smaller or equal energy).
In fact, the reverse is true as well:
$$\Cap_g (\phi^{-1} [A]) \leq C_d \Cap_g (A)$$
where the constant $C_d$ depends only on $d$, not on $g$ or $A$.
The proof of this fact is based on a trick of Benjamini
and Peres (1992b) and appears as Theorem~3.1 of
Pemantle and Peres (1995b).

Let $S$ be the random subtree of $4^\infty$ gotten from
the Bernoulli percolation with $p_n = n/(n+1)$, i.e.,
each vertex at level $n$ is killed with probability $1/(n+1)$.
The map $\phi$ carries $\partial S$ to a subset of
$[0,1]^2$ which is potential theoretically very similar 
to the range $G$ of Brownian motion run for a unit time, 
but which is easier to analyze because there is so
much independence.  This set may be used to derive 
properties of the intersections of independent Brownian motions.
Kakutani's Theorem says a single Brownian motion intersects
precisely those sets with positive logarithmic capacity,
so one may think of $G$ as having ``co-dimension'' log.  
This heuristic implies that the intersection of $k$ independent
copies of $G$ should have ``co-dimension'' $|\log|^k$,
which was conjectured in Tongring's thesis and first proved
by Fitzsimmons and Salisbury (1989).  The simplified proof
based on tree theory is due to Peres (1994) and goes as
follows.  

Let $A$ be any subtree of $4^\infty$.  The Sharp Bernoulli Theorem says 
the probability that $\partial A \cap \partial S$ is non-empty
is estimated up to constants by the capacity of $\partial A$
in a gauge $f (k) = 1/(k+1)$, which is equivalent in the sense of 
Section~2.3 to
the gauge $g(x) = |\log x|$.  Kakutani's theorem (or the quantitative
version found in Benjamini, Pemantle and Peres 1993) shows
that $\P (\phi [A] \cap G) \neq \emptyset)$ is estimated by 
$\Cap_{\log} (\phi [A])$, but this may be pulled back to 
$$\P (A \cap \phi^{-1} [G] \neq \emptyset ) \sim \Capb_f (A) .$$
Without abusing notation too much, we may identify subtrees of 
$4^\infty$ with closed subsets of $\partial 4^\infty$ and hence 
with closed subsets of $[0,1]^2$.
Thus $S$ is {\em intersection equivalent} to $G$ in the sense that
their probabilities of intersecting a third set differ by a bounded
factor.  Now several applications of Fubini's theorem finish
the proof as follows.  

Let $S_1$ and $S_2$ be IID copies of $S$ and $G_1$
and $G_2$ be IID copies of $G$.  The key fact is that while
$G' \deq G_1 \cap G_2$ is a mess, $S' \deq S_1 \cap S_2$ is a set with 
the same intersection estimates as $S$ except with $\log^2$ replacing 
$\log$.  Any set $A$ intersects $G'$ with positive probability
if and only if $A \cap G_1$ intersects $G_2$ with positive probability.
This is true if and only if $A \cap G_1$ intersects $S_2$ 
with positive probability, which is true if and only if $A \cap S_2$ 
intersects $G_1$ with positive probability.  This is true if and only if 
$A \cap S_2$ intersects $S_1$ with positive probability,
that is to say, if and only if $A$ intersects $S$ with positive 
probability, which we know to hold if and only if $\Cap_{\log^2} (A) > 0$.
Iteration extends this argument to common intersections with $k$ 
independent Brownian motions. 

Here is a similar argument due to Peres that proves Marstrand's
Intersection Theorem:

\vspace{1ex}\noindent{\bf Projection and Intersection Theorems 
(Marstrand 1954):~~}{\it (1) If $A \subseteq [0,1]^2$ is closed and
has positive 1-dimensional capacity, then $A$ intersects a 
random line with positive probability, where for specificity, 
we suppose the line is chosen by connecting two points chosen
independently and uniformly on the perimeter of $[0,1]^2$.
(2) For any $\ee > 0$, with positive probability the intersection 
has dimension at least $dim(A) - 1 - \ee$.}\\[1ex]
The proof of the Intersection Theorem~(2) from the Projection 
Theorem~(1) is as follows.  Suppose $dim(A) > 1$ and pick
any $\alpha \in (1 , dim(A))$.  The relationship between capacity
and dimension (Section~2.3) gives $\Cap_\alpha (A) > 0$.  
Claim: $\Cap_{\alpha - 1} (A \cap S) > 0$ with positive probability, 
where $S$ is a Galton-Watson subtree of $4^\infty$ corresponding to
percolation with probability $1/2$.  Proof of claim: Let
$S_\alpha$ and $S_{\alpha - 1}$ be independent Galton-Watson 
subtrees corresponding to percolation with probabilities $2^{-\alpha}$
and $2^{1 - \alpha}$ respectively.  Since $S \cap S_\alpha 
\diseq S_{\alpha - 1}$, the Sharp Bernoulli Theorem gives:
$\Cap_\alpha (A) > 0$ implies $\P (A \cap S_\alpha \neq \emptyset) > 0$
implies $\P ((A \cap S) \cap S_{\alpha - 1}) > 0$ which implies
that $\Cap_{\alpha - 1} (A \cap S) > 0$ with positive probability.
The claim is proved and~(2) follows immediately from the relation
between capacity and dimension.  For completeness we sketch a 
tree-based proof of the Projection Theorem (though Falconer (1990 page
103) considers the Intersection Theorem to be the more difficult of the two).
Identify $A$ with a subtree of $4^\infty$, let $l$ 
be the random subtree of $4^\infty$ corresponding to a random line, 
and notice that $W \deq l \cap A$ is quasi-Bernoulli.  The Basic
Second Moment Theorem then implies~(1).  

The correspondence $\phi$ also preserves Hausdorff measure
with respect to an arbitrary gauge.  This has been known
for at least 70 years, and was used by Frostman to prove
a lemma which is still of fundamental importance for the
analysis of fractal sets.  Frostman's lemma says that if
a subset $A$ of $[0,1]^d$ has positive Hausdorff measure 
with respect to a gauge $h$, then a measure $\mu$ exists 
for which $\mu (A) > 0$ and $\mu (C) \leq h (\mbox{diam}(C))$
for all sets $C$.  When transferred to the tree setting,
this becomes a special case of the min-cut-max-flow theorem.
All proofs I know (before 1994) of the 
existence of a Frostman measure translate the problem first to 
the tree setting (see for example Carleson 1967).  

The type of argument used by Frostman is very common, for example
see McMullen (1984), whose sequence space is transparently isomorphic 
to a regular tree.  The following example, mentioned in Section~1.2, 
illustrates the argument; see Falconer (1990) for another discussion of 
trees applied to random fractals.  

Let $C$ be a random set constructed from a measure $\mu$ on collections
of subcubes of $[0,1]^d$ as described in Section~1.2.  Let $\Gamma$ be
a tree representing $C$ as follows: the vertices of
$\Gamma_n$ are the chosen subcubes at the $n^{th}$ level of iteration,
and $X(v) = - \log r$ where $r < 1$ is the side of the subcube
corresponding to $v$ divided by the size of the subcube corresponding
to the parent of $v$.  The collection $\{ X(w) \}$ for all 
children $w$ of $v$ is IID as $v$ varies, which is enough independence
to apply results such as~(\ref{eq dim}).  It is easy to see that the 
dimension of $C$ is the same as the dimension of $\partial \Gamma$ in 
the metric $d(x,y) = \exp (-S(x \wedge y))$.  One may now proceed
directly via equation~(\ref{eq dim}) as follows.
One first computes the large deviation rate 
$$-g(\lambda) = \inf_a \{ -a \lambda + \log \E \exp (a X) \}$$  
and then applies~(\ref{eq dim}) to get 
$$\beta = \sup \{ \lambda - g(\lambda) : \lambda \leq a_0 \} .$$
An improvement is to use convex conjugate functions to 
see that the pair of optimizations leads to $\beta$ solving
$$ \E \sum X(v)^\beta = 1 , $$
summing over all children $v$ of a given vertex.  

In the deterministic case (self-similarity rather than stochastic 
self-similarity), $\sum X(v)^\alpha$ is constant and it well known 
that $\beta$ is the value of $\alpha$ that makes this equal to 1.  
Graf, Mauldin and Williams (1988) prove the stochastic version, 
giving as well the exact Hausdorff gauge function.  Their proof
of~(\ref{eq dim}) is long, but if one only wants the dimension of 
$\Gamma$ then Lyons' proof using percolation is best.

\section{Further direction, some theory, and some open problems}

The most important open problem about tree-indexed processes
is to get a converse to the Basic Theorem that would make
the quasi-Bernoulli case as sharp as the Bernoulli case.  

\begin{conj} \label{c1}
Given a tree-indexed process $\Omega, S, B, \mu$, let $f(n)
= \mu^n (\pi_n (B))^{-1}$ be the probability that a sequence of
$n$ IID picks from $\mu$ is extendable to some sequence in $B$.  
Then $P (\Gamma ; B) > 0$ implies $\Capb_f (\Gamma) > 0$.
\end{conj}

If this is true, then when $B$ is quasi-Bernoulli, this combines
with the Basic Theorem to show that $\Capb_f (\Gamma) > 0$ is
necessary and sufficient for $P (\Gamma ; B) > 0$.  The importance
in obtaining sharp results is that they can be used for the
sort of back-and-forth Fubini argument of the previous section.
One way to approach this conjecture is to try and understand the
nature of the event $A(\Gamma ; B)$ when it occurs.  For example,
does the existence of the witnessing path 
$(v_0 , v_1 , v_2 , \ldots )$ hinge on local luck 
(think of the existence of a point in a Poisson process 
of intensity one on $[0,1]$) or more a matter of there being 
plenty of chances globally (think of a supercritical branching process).  
The following fact is mentioned without proof in Pemantle and Peres (199?),
which proves the special case that $B$ is an increasing event.

\vspace{1ex}\noindent{\bf One-Implies-Many Theorem:~~}{\it
Let $\AB (\Gamma ; B)$ be the event that there exist uncountably
many paths $(v_0 , v_1 , v_2 , \ldots)$ for which 
$(X(v_1) , X(v_2) , \ldots) \in B$.  Suppose that $\mu^\infty (B) = 0$,
so that each fixed branch is a witness with probability zero.
Then $\P (A \setminus \AB) = 0$.  }\\[1ex]

In Fitzsimmons and Salisbury (1989) the necessity of the
capacity criterion is proved by showing that almost any
definition of first hitting time yields a measure with
finite energy.  Salisbury (1994, last page) asks for a similar
inequality for homogeneous trees. The following conjecture
is similar to Salisbury's problem, and would imply 
Conjecture~\ref{c1}:

\begin{conj} \label{c2}
Given a tree-indexed process, let $x(\omega)$ be any random 
element of $\partial \Gamma$ which is a witness to $A(\Gamma ; B)$
when $A$ occurs and is undefined otherwise.  Let $\nu$ be the
law of $x$, which is a subprobability measure of total mass
$P (\Gamma ; B)$.  If $P(\Gamma ; B) > 0$, then $\cEB_f (\nu) < \infty$,
with $f$ as in the previous conjecture.
\end{conj}

Another set of questions has to do with the domination relation
defined by $\Gamma \succeq \Gamma'$ if and only if for every
$\mu$ and $B$, $P (\Gamma ; B ; \mu) \geq P(\Gamma' ; B ; \mu)$.
This relation is understood at present only for spherically
symmetric trees and trees of height 2; see Pemantle and Peres (1994).
Understanding this even for trees of height 3 seems difficult. 
The notion of domination may be generalized to graded graphs as follows.
Say that $G$ is a graded graph if its vertex set may partitioned
into disjoint sets $V_0 , V_1 , \ldots , V_N$ such that edges
occur only between $V_j$ and $V_{j+1}$.  If $V_0 = \{ \rho \}$,
say that $\rho$ is the root of $G$.  For such a graded graph,
associate IID random variables $\{ X(v) \}$ to the vertices,
having common distribution $\mu$.
Let $P(G ; B)$ denote the probability that $(X(v_1) , \ldots ,
X(v_N)) \in B$ for some graded path in $G$ (a graded path being
a sequence of vertices $\{ v_j \}$ with $v_j \in V_j$ 
and consecutive vertices connected by edges).

\begin{conj} 
$P(G ; B) \leq P(G' ; B)$, where $G'$ is the graph consisting
of $M$ paths of length $N$ disjoint except at $\rho$ and 
$M$ is the number of distinct graded paths of $G$.  
\end{conj}

The Comparison Theorem implies this in the case where $G$ is
the graded graph of a tree.  An elementary proof in this
case is given by Benjamini and Peres (1992a).  Sidorenko 
(1991, 1992), motivated by the pursuit of other combinatorial
problems, proves special cases of this where $N=2$ 
and the graph is either acyclic or small.

Some tree-indexed processes that seem interesting in 
themselves are the tree-indexed Markov chains.  These 
were first studied by Benjamini and Peres (1992a, 1994b).  
Intuitively, these are branching Markov chains, where the
branching structure is prespecified as some fixed tree.  
To construct tree-indexed Markov chains
as standard tree-indexed processes, begin with transition
probabilities $p(x,y)$ on a countable space $Y$.
Let $\{ X(v) \}$ be IID uniform on the unit
interval, and for each $y \in Y$, let $\{ A_{y,z} : z \in Y \}$
be a partition of $[0,1]$ into sets such that the Lebesgue
measure of the set $A_{y,z}$ is $p(y,z)$.  Define $S(v)$ recursively
by $S(\rho) = y_0$ for some arbitrary $y_0$, and if $v$ is
the parent of $w$ then $S(w) = z$ if and only if $X(w) \in A_{S(v),z}$.
Thus along any single self-avoiding path from the root, one
sees a Markov chain with transitions $p (x,y)$.  

When the Markov chain is a random walk on a group, some
regularity of behavior can be established.
Benjamini and Peres (1994a) discuss the relation
between recurrence properties of such a tree-indexed random walk
and the growth of the group.  Amenability of $G$, for instance,
is equivalent to recurrence of the walk for any symmetric $\mu$
and any $\Gamma$ that has positive capacity in some gauge
$f(n) = e^{\alpha n}$, $\alpha > 0$ (these are just the
trees with positive Hausdorff dimension in the metric of Section~3.3).
Here, recurrence means the almost sure existence of infinitely
many $v$ with $S(v) = $ the identity.  Recurrence can be
determined from the Green's function if the group has polynomial
growth, but not if the group has a nontrivial Poisson boundary.

Many questions about tree-indexed random walks are still
open; here is just one.  Consider the case of a tree-indexed
random walk on $\R$.  Call $x = (v_0 , v_1 , \ldots)
\in \partial \Gamma$ an {\em escaping ray} if $S(v_n) \rightarrow
\infty$ as $n \rightarrow \infty$.  Say that $x$ is a {\em bouncing ray}
if $\infty > \liminf S(v_n) > -\infty$.  

\begin{conj} \label{c4}
If $\mu$ has mean zero and finite variance and
$\Gamma$ almost surely has bouncing rays, then $\Gamma$
almost surely has escaping rays.
\end{conj}

This is known to be true only in the cases where the $\mu$
is normal or Rademacher (Pemantle and Peres 1995a).  The
proofs in these cases involve a long detour through 
potential-theoretic equivalences, conspicuously absent 
in the statements.  Probably the conjecture is true for the
reason that $\Capb_f (\Gamma) > 0$ is necessary and sufficient
in both cases, where $f (k) = \sqrt{k}$.  Is there an 
elementary argument?  

Finally, as I promised in Section~1.3, I will sketch a connection
between random walks on trees and the potential theory of trees.
A very detailed such connection may be made via electrical network
theory, but the connection via potential theory seems more germane here.

Suppose we wish to estimate the probability that a simple random walk 
started from the root of a finite tree $\Gamma$ of height $N$, hits
the boundary of the tree before hitting a cemetery $\Delta$ attached
to the root.  The set $W$ of vertices of $\Gamma_N$ hit by a random
walk before hitting $\Delta$ is a random set satisfying~(\ref{eq QB}) 
of the Basic Second Moment Theorem.  It also satisfies a certain
Markov property, and these two facts imply that this random set is
intersection equivalent to a Bernoulli percolation; the cumulative 
survival probabilities for this percolation turn out to be $1/(n+1)$ 
at level $n$, corresponding to a kernel $K(x,y) = |x \wedge y|$.
Letting the height $N \rightarrow \infty$ gives a criterion 
for transience discovered by Lyones (1990) 
(see also Benjamini and Peres 1992b):
simple random walk on $\Gamma$ is transient if and only if
$\partial \Gamma$ has positive $|x \wedge y|$-capacity.   
 
\noindent{\sc Acknowledgements:}  This paper would not have been
written were it not for David Griffeath, would not have been 
suitable for this journal were it not for Rob Kass, and would
not have been any good were it not for Yuval Peres.  

\renewcommand{\baselinestretch}{1.0}\large 
\vspace{.4in}

\today

\end{document}